\documentclass{amsart}
\usepackage{amssymb}
\usepackage{eucal}
\usepackage{amsfonts}
\usepackage{epsfig}
\usepackage{hyperref}
\usepackage{color}
\usepackage[all]{xy}
\let\oldcite\cite                                  
\newtheorem{thm}{Theorem}[section]
\newtheorem{cor}[thm]{Corollary}
\newtheorem{lem}[thm]{Lemma}
\newtheorem{prop}[thm]{Proposition}
\theoremstyle{definition}
\newtheorem{defn}[thm]{Definition}
\theoremstyle{remark}
\newtheorem{rem}[thm]{Remark}
\numberwithin{equation}{section} \theoremstyle{remark}

\newcommand{\X}{\mathcal{X}}
\newcommand{\Y}{\mathcal{Y}}
\newcommand{\Z}{\mathcal{Z}}
\newcommand{\K}{\mathcal{K}}

\newcommand{\bbX}{\mathbb{X}}
\newcommand{\bbY}{\mathbb{Y}}

\let\:=\colon

\newcommand{\Top}{\mathsf{Top}}
\newcommand{\CGTop}{\mathsf{CGTop}}
\newcommand{\TopSt}{\mathfrak{TopSt}}

\newcommand{\St}{\mathfrak{St}}

\newcommand{\ParSt}{\mathfrak{ParSt}}

\newcommand{\tX}{\tilde{\X}}

\newcommand{\Hom}{\operatorname{Hom}}
\newcommand{\Map}{\operatorname{Map}}

\newcommand{\id}{\operatorname{id}}
\newcommand{\pr}{\operatorname{pr}}
\newcommand{\ev}{\operatorname{ev}}
\newcommand{\hFib}{\operatorname{hFib}}

\def\smashedst{\setbox0=\hbox{$\rightrightarrows$}\ht0=4pt\box0}
\newcommand{\sst}[1]{\stackrel{#1}{\smashedst}}

\begin{document}

\title{Mapping stacks of topological stacks}

\author{Behrang Noohi}

\begin{abstract}
 We prove that the mapping stack $\Map(\Y,\X)$ of topological
 stacks $\X$ and $\Y$ is again a topological stack if $\Y$ admits a
 groupoid presentation $[Y_1\sst{}Y_0]$ such that $Y_0$ and $Y_1$
 are compact topological spaces. If $Y_0$ and $Y_1$ are only
 locally compact, we show that $\Map(\Y,\X)$ is a paratopological stack.
 In particular, it has a classifying space (hence, a natural weak 
 homotopy type).
 We also show that the weak homotopy type of the  mapping stack
 $\Map(Y,\X)$ does not change if we replace $\X$ by its classifying space,
 provided that $Y$ is a paracompact topological space.  As an example, we 
 describe the loop stack  of the classifying stack $\mathcal{B}G$ of a topological 
 group $G$ in terms of twisted loop groups of $G$.
\end{abstract}

\maketitle%

\tableofcontents%

\section{Introduction}

Let $\X$ and $\Y$ be topological stacks. The purpose of these notes is
to show that under a mild condition on $\Y$ the mapping stack
$\Map(\Y,\X)$ is a paratopological stack and, in particular, admits a
classifying space.

There are various classes of stacks to which this result applies. For example, 
for arbitrary $\X$, 
we can take $\Y$ to be coming from a:
Lie groupoid, orbifold,
action of a locally compact group on a locally compact space (e.g., 
classifying stack $\mathcal{B}G$ of locally compact
group),  complex-of-groups,  Artin stack of finite type over
complex numbers,  foliation on a manifold, and so on.

In the case where $\X$ and $\Y$ are orbifolds, the mapping stack
$\Map(\Y,\X)$ has been studied by Chen. One of the main results of
\cite{Chen} is that in this case $\Map(\Y,\X)$ is again an orbifold.
To our knowledge, this is the only general result previously known 
about the mapping
stacks being  topological, and its proof is quite nontrivial.
Another known case is when $\Y=S^1$, in which case the free loop
stack $\Map(S^1,\X)=:\mathcal{L}\X$ is  shown to be a topological
stack in \cite{LuUr}.

The mapping stack $\Map(\Y,\X)$ can be defined for arbitrary stacks
$\X$ and $\Y$, and it is functorial in both $\X$ and $\Y$
($\S$\;\ref{S:Generalities}). However, it does not in general admit a groupoid
presentation, even if $\X$ and $\Y$ do. Therefore,  $\Map(\Y,\X)$ may 
not always be 
a topological stack. 

In the case where $\X$ and $\Y$ are
topological spaces, $\Map(\Y,\X)$ coincides with the usual mapping
space with the compact-open topology. If we are given groupoid
presentations $\bbX=[X_1\sst{}X_0]$ and $\bbY=[X_1\sst{}X_0]$ for
$\X$ and $\Y$, the mapping stack $\Map(\Y,\X)$ parameterizes the
Hilsum-Skandalis morphisms from $\bbY$ to $\bbX$. In the case
$\X=\mathcal{B}G$, where $G$ is an arbitrary topological group, the
mapping stack $\Map(\Y,\mathcal{B}G)$ classifies principal
$G$-bundles over $\Y$.

Our first main result shows that, under a locally compactness condition on $\Y$,
the mapping stack $\Map(\Y,\X)$ admits a classifying space 
(Definition \ref{D:classifying}),
hence also a natural weak homotopy type; see 
Theorems \ref{T:mapping1} and \ref{T:mapping2}.

\begin{thm}{\label{T:1}}
  Let $\X$ and $\Y$ be topological stacks, and let $\Map(\Y,\X)$ be
  their mapping stacks. If $\Y$ admits a presentation by a groupoid
  $[Y_1\sst{}Y_0]$ such that $Y_0$ and $Y_1$ are compact topological
  spaces, then $\Map(\Y,\X)$ is again a topological stack. If
  $Y_0$ and $Y_1$ are only locally compact, then $\Map(\Y,\X)$ is
  paratopological (Definition \ref{D:paratop}).
\end{thm}

The theorem implies that
the mapping stack $\Map(\Y,\X)$ has a {\em classifying space} 
(Definition \ref{D:classifying}). In other words, there exists a topological
space $V$ and a map $\varphi \: V \to \Map(\Y,\X)$ that is a {\em universal
weak equivalence}, i.e., $\varphi$ has the property that
for every map $T  \to \Map(\Y,\X)$ from 
a topological space $T$, 
the base extension
$\varphi_T \: V_T \to T$ of $\varphi$ is a weak equivalence of
topological spaces. (In fact, whenever $T$ is paracompact, we can
arrange so that $\varphi_T$ admits a section and $V_T$ has a fiberwise 
strong
deformation retraction over the image of this section.) In the case
where we have the compactness condition on $\Y$, we can arrange so
that $\varphi$ is also an epimorphism (thereby, giving rise to a
presentation for $\Map(\Y,\X)$ by the topological groupoid
$[U\sst{}V]$, where $U=V\times_{\Map(\Y,\X)}V$).

As is discussed in \cite{Homotopytype} in detail, existence of such a 
classifying space
$\varphi \: X \to \X$ is crucial for doing algebraic topology on
$\X$, as it allows one to translate problems about the stack $\X$ to ones
about the space $X$ (e.g., by pull-back along $\varphi$).  
For example, classifying spaces of stacks have been used 
extensively in \cite{BGNX} to develop 
intersection theory on loop stacks of differentiable stacks. Another
application appears in \cite{EbGi} in which the same method is used
to produce new classes in the singular homology (with coefficients)
of moduli stacks $\bar{\mathfrak{M}}_{g,n}$ of curves.

As the terminology suggests, in the case where $\X$ is the quotient stack of a 
topological groupoid $\bbX=[X_1\sst{}X_0]$, the Haefliger classifying space
of $\bbX$ is indeed a classifying space for $\X$ in the sense discussed above.
Therefore, the weak homotopy type of $\X$ is the same as the weak homotopy type 
of the Haefliger classifying space.
This raises the question of {\em homotopy invariance} of mapping stacks. For example,
one could ask whether the loops stack $L\X$ has the same weak homotopy type
as the loop space $LX$ of the classifying
space of $\X$. In $\S$\ref{S:Invariance} we give an affirmative answer 
to this question. More generally, we prove the following 
(see Corollary \ref{C:invariance1}).

\begin{thm}
  Let $\X$ be  a topological stack, and let $X$ be a classifying space for it. Let $Y$ be a 
  paracompact topological space. Then, there is a natural (universal) weak equivalence 
  $\Map(Y,X) \to \Map(Y,\X)$. That is, $\Map(Y,X) $ is a classifying space for $\Map(Y,\X)$.
\end{thm}

As another application of Theorem \ref{T:1}, we see in $\S$\ref{SS:htpyfiber} 
that a morphism $f \: \X \to \Y$  of topological stacks  factorizes as a composition 
$f=p_f\circ i_f$ of a closed embedding $i_f \: \X \to \tilde{\X}$ which admits a 
strong deformation retraction followed by a Hurewicz fibration $p_f \: \tilde{\X} \to \Y$;
see Proposition \ref{P:replacement}. In particular, one can define the homotopy fiber
of $f$ as a topological stack.

Finally, we study the loop stack of the classifying stack 
$\mathcal{B}G$ of a topological group $G$.  We prove the following.

\begin{thm} 
   Let $G$ be a topological group (not necessarily connected). 
   Then, there are natural  weak homotopy  equivalences
     $$LBG\cong L\mathcal{B}G \cong \coprod_{i\in C_G} BL_{(\alpha_i)}G.$$ 
   Here, $BG$ is the Milnor classifying space of $G$, $C_G$ is the set of conjugacy 
   classes of $\pi_0G$, and  $L_{(\alpha_i)}G$
   are twisted loop groups.
   In particular, when $G$ is discrete, there are natural weak homotopy equivalences
     $$LBG\cong \mathcal{IB}G \cong G\times_GEG,$$ 
   where  $\mathcal{IB}G$ is the inertia stack and $G\times_GEG$ is
   the Borel construction on the conjugation action of $G$ on itself. (The left equivalence
   is indeed an equivalence of stacks.)
\end{thm}

There is also a similar description for the weak homotopy type of the pointed
loop stack $\Omega\mathcal{B}G$; see Theorem \ref{T:pointed}.

\bigskip
\noindent{\bf Conventions.} Throughout the notes, $\CGTop$ stands
for the category of compactly generated topological spaces. All
topological spaces will be assumed to be compactly generated. All
stacks considered are over $\CGTop$.

\section{Topological and paratopological stacks}{\label{S:Topological}}
 
In this section, we  recall some  basic facts and definitions from  
\cite{Foundations} and
\cite{Homotopytype}. 

By a {\bf topological stack} we mean a stack $\X$ over $\CGTop$ which is
equivalent to the quotient stack of a topological groupoid
$\bbX=[X_1\sst{}X_0]$ with $X_1$ and $X_0$ (compactly generated)
topological spaces.

There are two classes of stacks that are of special interest to us in this paper: 
paratopological stacks (that are more general than topological stacks) and 
Hurewicz topological stacks (that are special types of topological stacks).
We will discuss them shortly.

\subsection{Classifying space of a stack}{\label{SS:classifying}}

The following notion plays an important role in the homotopy theory of topological stacks.

\begin{defn}[\cite{Homotopytype}, Definition 5.1]{\label{D:univwe}}
 A representable morphism $f \: \X \to \Y$ of stacks is called a {\bf
 universal weak equivalence} if the base extension
 $f_T \: \X_T \to T$ of $f$ to an arbitrary topological space $T$ is
 a weak equivalence of topological spaces.
\end{defn}

\begin{defn}{\label{D:classifying}}
  Let $\X$ be a stack whose diagonal $\X \to \X\times\X$ is representable. 
  By a {\bf classifying space} for $\X$ we mean a topological space $X$ 
  equipped with a
  universal weak equivalence $\varphi \: X \to \X$.
\end{defn}

In \cite{Homotopytype}, a stack $\X$ which admits a classifying space is called a 
{\em homotopical stack}.  The classifying space of a homotopical stack 
$\X$ is unique up to a 
unique isomorphism
in the weak homotopy category of topological spaces. It is functorial (in 
the weak homotopy category) and is a model for the weak homotopy type 
of $\X$.

\subsection{Paratopological stacks}{\label{SS:paratopological}}

\begin{defn}[\oldcite{Homotopytype}, Definition 9.1]{\label{D:paratop}}
 We say that a stack $\X$ is {\bf paratopological} if it
 satisfies the following conditions:
 \begin{itemize}
   \item[$\mathbf{A1.}$] Every map $T \to \X$ from a topological
     space $T$ is representable (equivalently, the diagonal
     $\X \to \X\times\X$ is representable);

   \item[$\mathbf{A2.}$] There exists a
      morphism $X \to \X$ from a topological space $X$ such that for
      every morphism $T \to \X$, with $T$ a paracompact topological
      space, the base extension $T\times_{\X} X
       \to T$ is an epimorphism (i.e., admits local sections).
 \end{itemize}
\end{defn}

We denote the 2-categories of stacks, topological stacks, and paratopological
stacks by $\St$, $\TopSt$, and $\ParSt$, respectively. We have  full inclusions of 
2-categories 
       $$\TopSt \subset \ParSt \subset \St.$$

The following lemma says that a paratopological stacks looks like a 
topological stack in the eye of a paracompact topological space.

\begin{lem}[\oldcite{Homotopytype}, Lemma 9.2]{\label{L:approximation}}
 Let $\X$ be a stack such that the diagonal $\X \to
 \X\times\X$ is representable. Then, $\X$ is paratopological
  if and only if there exists a
 topological stack $\bar{\X}$ and a representable morphism $p \: \bar{\X} \to \X$
 such that for every paracompact topological space $T$, $p$ induces an equivalence of
 groupoids  $\bar{\X}(T) \to \X(T)$.
\end{lem}

\begin{proof}
  We just indicate how $\bar{\X}$ is constructed. Let $X \to \X$ be as in 
  Definition \ref{D:paratop} ($\mathbf{A2}$). Set $X_0:=X$ and $X_1:=X\times_{\X}X$.
  The quotient stack $\bar{\X}$ of the groupoid $[X_1\sst{}X_0]$ has the desired
  property. 
\end{proof}

  Note that paracompactness of $T$ does not play a role in the above proof; if in
  Definition \ref{D:paratop} ($\mathbf{A2}$) the space $T$ is required to belong to a certain
  class of spaces (e.g., paracompact, CW complex, etc.), then Lemma \ref{L:approximation} 
  will be true with $T$ in the same class of spaces.

\begin{defn}[\oldcite{Homotopytype}, Definition 5.1]{\label{D:shrinkable}}
A representable morphism $f \: \X \to \Y$ of stacks is called {\bf
parashrinkable} if for every morphism $T \to \Y$ from a paracompact
topological space $T$, 
the base extension $f_T \: \X_T \to T$ of $f$ over $T$ admits a
section $s$ and a fiberwise deformation retraction of $\X_T$ onto
$s(T)$.
\end{defn}

 \begin{prop}[\oldcite{Homotopytype}, Proposition 9.4]{\label{P:classifying}}
  For every paratopological stack $\X$, there exists a parashrinkable
  morphism $\varphi \: X \to \X$ from a topological space $X$. If $\X$
  is a topological stack, then such a $\varphi$ can be chosen to be an
  epimorphism.
 \end{prop}

Observe that a parashrinkable morphism  is a 
universal weak equivalence (Definition \ref{D:univwe}). Therefore,
the space $X$ in the above proposition is a classifying space 
for $\X$ in the sense
of Definition \ref{D:classifying}. Also, note that every topological
stack $\X$ is a paratopological
stack, hence admits a classifying space by the above proposition.  
In fact, it is shown in
(\cite{Homotopytype}, Theorem 6.1) that the Haefliger classifying space of 
a groupoid presentation
for $\X$ is a classifying space for $\X$ in the sense of Definition \ref{D:classifying}.
Let us record this fact for future reference.

\begin{prop}[\oldcite{Homotopytype}, Theorem 6.1]{\label{P:Haefliger}}
  Let $\X$ be a topological stack and  $\bbX=[X_1\sst{}X_0]$ a groupoid presentation
  for it. Let $B\bbX$ denote the Haefliger classifying space of $\bbX$. Then there
  is a natural parashrinkable morphism 
       $$\varphi \: B\bbX \to \X.$$
  In particular, for every topological group $G$, there is a parashrinkable
  morphism 
      $$\varphi \: BG \to \mathcal{B}G,$$
 where $BG$ is the Milnor classifying space of $G$. 
\end{prop}

\subsection{Hurewicz topological stacks}{\label{SS:Hurewicz}}

In the 2-category of topological stacks colimits are quite ill-behaved. 
(For example, homotopies between maps with 
target $\X$ a topological stack do not always glue.)
The situation is slightly better in the 2-category of
 {\em Hurewicz topological stacks} (e.g., see Proposition \ref{P:glue}). 
This makes Hurewicz topological stacks more appropriate for certain 
homotopy theoretic manipulations.

To give the definition of a Hurewicz topological stack, we recall 
some standard definitions. 
A {\em Hurewicz fibration} is a continuous map of topological spaces
which has the homotopy lifting
property for all topological spaces. A map $f \: X \to Y$ of
topological spaces is a {\em local Hurewicz fibration} if for every
$x \in X$ there are opens $x \in U$ and $f(x) \in V$ such that
$f(U)\subseteq V$ and $f|_U\: U \to V$ is a Hurewicz fibration.  An
important example is that  of a topological
submersion: a map $f:X\to Y$, such that locally $U$ is homeomorphic
to $V\times \mathbb{R}^n$, for some $n$.

\begin{defn}\label{D:Hurewicz}
   A topological stack $\X$ is called {\bf Hurewicz} if it
   is equivalent to the quotient stack $[X_0/X_1]$ of a topological
   groupoid $[X_1\sst{} X_0]$ whose source and target maps are local
   Hurewicz fibrations.
\end{defn}

Hurewicz topological stacks form a full sub 2-category of $\TopSt$.

\subsection{Limits of topological stacks}{\label{SS:limits}}

We will need the following fact in the proof of Theorem
\ref{T:mapping2}.

\begin{prop}{\label{P:arblimits}}
  The 2-categories $\TopSt$, $\ParSt$, and the 2-category of 
  Hurewicz topological stacks are closed under finite limits.
  The 2-category of stacks with representable diagonal is closed under arbitrary limits.
  The 2-category $\ParSt$ is closed under arbitrary fiber products.
  In particular, the product of an arbitrary family of
  paratopological stacks is paratopological.
\end{prop}

\begin{proof}
  The statements about $\TopSt$ and $\ParSt$ are proved in 
  \cite{Homotopytype}, Lemmas 9.12 and 9.13. The same argument used in 
  [ibid.] proves that Hurewicz topological stacks are closed under finite limits.
  The statement about stacks with representable diagonal is 
  (\cite{Homotopytype}, Lemmas 9.11).
\end{proof}

\section{Generalities on mapping stacks}{\label{S:Generalities}}

We begin by recalling the definition of the mapping stack. Let $\X$
and $\Y$ be   stacks over $\CGTop$. We define the stack
$\Map(\Y,\X)$, called the {\bf mapping stack} from $\Y$ to $\X$, by
the rule
  $$ T\in \CGTop \ \ \   \mapsto  \ \ \ \Hom(T\times \Y,\X)\,,$$
where $\Hom$ denotes the groupoid of stack morphisms. This is easily
seen to be a stack. 
Note that we have a natural equivalence of groupoids
    $$\Map(\Y,\X)(*)\cong\Hom(\Y,\X),$$
where $*$ is a point. In particular, the underlying set of the coarse 
moduli space of 
$\Map(\Y,\X)$ is the set of 2-isomorphism classes of morphisms
from $\Y$ to $\X$.

It follows from the exponential law for mapping
spaces that when $X$ and $Y$ are spaces, then
$\Map(Y,X)$ is representable by the usual mapping space from $Y$ to
$X$ (endowed with the compact-open topology).

The mapping stacks are functorial in both variables.

\begin{lem}{\label{L:functorial}}
 The mapping stacks $\Map(\Y,\X)$ are functorial in $\X$ and
 $\Y$.
 That is, we have natural functors
 $\Map(\Y,-) \: \St \to \St$ and $\Map(-,\X) \: \St^{op} \to
 \St$. (Note: in $\St^{op}$ we only invert the direction of
 1-morphisms.)
\end{lem}

The exponential law holds for mapping stacks.

\begin{lem}{\label{L:exponential}}
 For stacks $\X$, $\Y$ and $\Z$ we have a natural equivalence of stacks
   $$\Map(\Z\times\Y,\X)\cong\Map(\Z,\Map(\Y,\X)).$$
\end{lem}

\section{Mapping stacks of topological stacks}{\label{S:Mapping}}

In this section we show that, for a fairly large class of
topological stacks $\X$ and $\Y$, the machinery developed in
\cite{Homotopytype} can be used to associate a classifying space 
(hence, a homotopy type) to the
mapping stack $\Map(\Y,\X)$.

The first main result of this section (Theorem \ref{T:mapping1})
shows that if $\X$ and $\Y$ are topological stacks, then $\Map(\Y,\X)$ is 
again a topological stack, provided that $\Y$ has a
groupoid presentation $[Y_1\sst{}Y_0]$ in which both $Y_1$ and $Y_0$
are compact topological spaces. The compactness assumption is
somewhat restrictive, although it is enough for many applications (e.g., loop
stacks). In Theorem \ref{T:mapping2} we show that
$\Map(\Y,\X)$ is a {\em paratopological} stack under the weaker
assumption that $Y_0$ and $Y_1$ are only locally compact. By the results
of \cite{Homotopytype},  this allows one to associate a classifying space 
to $\Map(\Y,\X)$.

\begin{lem}{\label{L:A1}}
   Let $\X$ and $\Y$ be  stacks. Assume that $\Y$ is topological
   and that the diagonal $\X \to \X\times\X$ is representable. Then,
   for every topological space $S$, every morphism $S \to
   \Map(\Y,\X)$ is representable. (Equivalently,  
   $\Map(\Y,\X)$ has a representable diagonal.)
\end{lem}

\begin{proof}
 First note that we can reduce to the case where $\Y=Y$ is a
 topological space. This is possible because by (\cite{Foundations},
 Proposition 3.19) the mapping stack $\Map([Y_0/Y_1],\X)$ can be
 written as the limit of a (finite) diagram produced out of
 $\Map(Y_0,\X)$ and $\Map(Y_1,\X)$; see Proposition \ref{P:arblimits}.

 Let $p \: S \to \Map(Y,\X)$ and $q \: T \to \Map(Y,\X)$ be
 arbitrary morphisms from  topological spaces $S$ and $T$. We
 need to show that $T\times_{\Map(Y,\X)}S$ is a topological space.
 Let
 $\tilde{p} \: S\times Y \to \X$ be the defining map for $p$, and
 $\tilde{q} \: T\times Y \to \X$ the one for $q$. Set $Z:=(T\times
 Y)\times_{\X}(S\times Y)$. This is a topological space
 which sits in the following 2-cartesian diagram:
   $$\xymatrix{ \Map(Y,Z) \ar[r] \ar[d]
                     & \Map(Y,S\times Y)  \ar[d]^{\tilde{p}_*}\\
        \Map(Y,T\times Y) \ar[r]^{\tilde{q}_*} & \Map(Y,\X)    }$$
 The claim now follows from the fact that $p$ and $q$ factor through
 $\tilde{p}_*$ and $\tilde{q}_*$, respectively.
\end{proof}

\begin{thm}{\label{T:mapping1}}
  Let $\X$ and $\K$ be  topological stacks. 
  Assume that
  $\K\cong[K_0/K_1]$, where $[K_1\sst{}K_0]$ is a topological
  groupoid with $K_0$ and $K_1$ compact.
  Then, $\Map(\K,\X)$ is a topological stack.
\end{thm}

\begin{proof}
 First note that we can reduce to the case where $\K=K$ is a compact
 topological space. This is possible because by (\cite{Foundations}, Proposition
 3.19) the mapping stack $\Map([K_0/K_1],\X)$ can be written as the
 limit of a finite diagram produced out of $\Map(K_0,\X)$ and
 $\Map(K_1,\X)$; use Proposition \ref{P:arblimits}.

 In view of Lemma \ref{L:A1}, all we need to do is to find an epimorphism 
 $R \to \Map(K,\X)$ from a topological space $R$. First some
 notation. Let $\mathbb{Y}=[Y_1\sst{}Y_0]$ and
 $\mathbb{X}=[X_1\sst{}X_0]$ be topological groupoids. 
 We define $\Hom(\mathbb{Y},\mathbb{X})$ to be the
 space of continuous groupoid morphisms from $\mathbb{Y}$ to
 $\mathbb{X}$. This is topologized as a subspace of
 $\Map(Y_1,X_1)\times\Map(Y_0,X_0)$, and
 it represents the (set-valued) functor
  $$T\in \Top \ \ \   \mapsto  \ \ \ \text{groupoid morphisms}
              \ \ T\times \mathbb{Y} \to \mathbb{X},$$
 where $T\times \mathbb{Y}$ stands for the groupoid $[T\times
 Y_1\sst{}T\times Y_0]$. In particular, we have a universal family of
 groupoid morphisms
 $\Hom(\mathbb{Y},\mathbb{X})\times\mathbb{Y} \to \mathbb{X}$.

 Let  $J$ be the set of all finite open covers of $K$. For $\alpha
 \in J$, let $U_{\alpha}$ denote the disjoint union of the open sets
 appearing in the open cover  $\alpha$. There is a natural  map
 $U_{\alpha} \to K$. Let $\mathbb{K}_{\alpha}:=[U_{\alpha}\times_K
 U_{\alpha}\sst{}U_{\alpha}]$ be the corresponding topological
 groupoid. 
 Note that the quotient stack of $\mathbb{K}_{\alpha}$ is $K$.
 Fix a  groupoid presentation $\mathbb{X}=[X_1 \sst{} X_0]$   for
 $\X$, and let $\pi \: X_0 \to \X$ be the corresponding atlas for
 $\X$. Set $R_{\alpha}=\Hom(\mathbb{K}_{\alpha},\mathbb{X})$, 
 with $\Hom$ being  as above. Let
 $R=\coprod_{\alpha} R_{\alpha}$.

 The universal groupoid morphisms $R_{\alpha} \times
 \mathbb{K}_{\alpha} \to \mathbb{X}$ give rise to morphisms
 $R_{\alpha} \to \Map(K,\X)$. Putting these all together we obtain a
 morphism $R \to \Map(K,\X)$. We claim that this an epimorphism.
 (Here is where compactness of $K$ gets used.) Let $p \: T \to
 \Map(K,\X)$ be  an arbitrary morphism. We have to show that, for
 every $t \in T$, there exists  an open neighborhood $W$ of $t$ such
 that $p|_W$ lifts to $R$. Let $\tilde{p} \: T\times K \to \X$ be
 the defining morphism for $p$. Since $\pi: X_0 \to \X$ is an
 epimorphism, we can find finitely many open sets $V_i$ of $T\times
 K$ which cover $\{t\}\times K$ and such that $\tilde{p}|_{V_i}$
 lifts to $X_0$ for every $i$. We may assume $V_i=K_i\times W$,
 where $K_i$ are open subsets of $K$, and $W$ is an open
 neighborhood of $t$ independent of $i$. Let $\alpha:=\{K_i\}$  be
 the corresponding open cover of $K$. Then $p|_W$ lifts
 to $R_{\alpha} \subset R$.
\end{proof}

\begin{rem}
  In the above proof we implicitly made use of the fact that the
  cartesian product of a compactly generated topological space with
  a compact topological space is again compactly generated
  (\cite{Whitehead}, I.4.14.). (Recall that product in
  $\CGTop$ is, in general, slightly different from the usual
  cartesian product in $\Top$. This is due to the fact that the product of two
  compactly generated spaces may no longer be compactly generated.
  See (\cite{Whitehead}, I.4) for more details.)
\end{rem}

We now treat the case of mapping stacks $\Map(\Y,\X)$ where $\Y$ is
no longer compact. We have the following result.

\begin{thm}{\label{T:mapping2}}
  Let $\X$ and $\Y$ be paratopological stacks. Assume that
  $\Y\cong[Y_0/Y_1]$, where $[Y_1\sst{}Y_0]$ is a topological
  groupoid with $Y_0$ and $Y_1$ locally compact.
  Then, $\Map(\Y,\X)$ is a paratopological
  stack (Definition \ref{D:paratop}). In particular,
  $\Map(\Y,\X)$ has a classifying space  in the sense of Definition \ref{D:classifying}.
\end{thm}

\begin{proof}
   We prove the theorem in several steps.
   
   \medskip \noindent {\em Step 1, $\X$ a paratopological stack, and 
   $\Y=Y$ a compact topological space.}
    If $\X$ is topological, then $\Map(Y,\X)$ is topological by  
   Theorem \ref{T:mapping1}. For an arbitrary paratopological stack $\X$, pick 
   a topological stack $\bar{\X}$ as in Lemma \ref{L:approximation}. For 
   every paracompact space $T$, the product $T\times Y$ is paracompact. This 
   implies that, $\Map(Y,\bar{\X})(T) \to  \Map(Y,\X)(T)$ is an equivalence of groupoids
   for every paracompact $T$. Furthermore, by Lemma \ref{L:A1}, $\Map(Y,\X)$ has 
   representable diagonal. It follows from Lemma \ref{L:approximation} that
   $\Map(Y,\X)$ is a paratopological  stack.   
   
   \medskip \noindent {\em Step 2, $\X$ a paratopological stack, and $Y$ a 
   disjoint union of compact topological spaces.}  Write $Y=\coprod Y_i$, with $Y_i$ compact. 
   It is  easy to see that $\Map(\coprod Y_i,\X)\cong\prod\Map(Y_i,\X)$. Hence, by 
   Proposition \ref{P:arblimits}, $\Map(Y,\X)$ is a paratopological stack.

  \medskip \noindent {\em Step 3, $\X$ a paratopological stack, and $\Y=[Y_0/Y_1]$ 
   quotient stack of 
   a topological groupoid $[Y_1 \sst{} Y_0]$ such that $Y_0$ and $Y_1$ are 
   disjoint unions of compact topological spaces.}
   As in the proof of Theorem \ref{T:mapping1},
   the mapping stack $\Map([Y_0/Y_1],\X)$ can be written as the
   limit of a finite diagram produced out of $\Map(Y_0,\X)$ and
   $\Map(Y_1,\X)$. Now use Proposition \ref{P:arblimits} and Step (2).
  
  \medskip \noindent {\em Step 4, $\X$ a paratopological stack, and $\Y=Y$ a 
   locally compact topological space.}  Choose a collection $\{Y_i\}$
   of closed compact subsets of $Y$ whose interiors covers $\Y$, and let $Y_0$ 
   denote their disjoint
   union. Let $Y_1:=Y_0\times_{Y}Y_0$ be the disjoint union of pairwise intersections of 
   the compact sets $Y_i$. Then $Y$ is the quotient stack
   of the groupoid $[Y_1 \sst{} Y_0]$ and the claim follows from Step (3).

  \medskip \noindent {\em Step 5, $\X$ a paratopological stack, and $\Y=[Y_0/Y_1]$ 
   quotient stack of 
   a topological groupoid $[Y_1 \sst{} Y_0]$ such that $Y_0$ and $Y_1$ are locally 
   compact topological spaces.} Use the same argument as in Step (3).
\end{proof}

\begin{cor}{\label{C:diff}}
  Let $\Y$ be a differentiable stack and $\X$ and arbitrary
  paratopological stack. Then, $\Map(\Y,\X)$ is a paratopological
  stack. In particular, $\Map(\Y,\X)$ admits a classifying space.
\end{cor}

\section{Application}{\label{S:Application}}

In this section we present an application of our results by showing that 
any morphism 
$f \: \X \to \Y$ of topological stacks has a factorization as a homotopy 
equivalence followed by a Hurewicz fibration. Along the way, we also 
point out certain subtleties that one should be aware of when working 
with mapping stacks. This explains why in certain contexts 
it is preferable to work with
Hurewicz topological stacks as opposed to arbitrary topological stacks.

\subsection{A gluing lemma for mapping stacks}{\label{SS:gluing}}

One main subtlety of working with mapping stacks is that certain intuitive 
statements about them may not hold in full generality.
The following proposition is proved in (\cite{BGNX}, Proposition 2.2). 
The proof is easy and it 
relies on Proposition 1.3 of [ibid.] (see \cite{Foundations}, Theorem 16.2 for an 
earlier version of 
this proposition). We emphasize that the statement may not be true if we do not
assume that $\X$ is a {\em Hurewicz} topological stack.

\begin{prop}{\label{P:glue}}
  Let $j \: A\to Y$ be a closed embedding of Hausdorff spaces. Assume that $j$ is a
 local cofibration in the sense that every $x \in A$ has a neighborhood $U$ in $Y$ such that
 $j|_{A \cap U} \: A \cap U \to U$ is a Hurewicz cofibration.
 Let $A\to Z$ be a finite proper map of Hausdorff spaces.
  Let $\X$ be a Hurewicz topological stack. Then the  diagram
     $$\xymatrix{
       \Map(Z\vee_A Y, \X) \ar[r] \ar[d] &  \Map(Y, \X) \ar[d] \\
       \Map(Z, \X) \ar[r]   & \Map(A, \X)   }$$
  is a 2-cartesian diagram  of stacks.
\end{prop}

\subsection{Path and loop stacks}{\label{SS:path}}

It is immediate from Theorem \ref{T:mapping1} that the {\bf path stack}
$P\X=\Map([0,1],\X)$ of a topological (resp., paratopological) stack $\X$ is again
a topological (resp., paratopological) stack. There is a
 natural {\em constant path} morphism $c_{\X} \: \X \to P\X$. For
 for every $t \in [0,1]$, there is a
 natural evaluation map $\ev_t \: P\X \to \X$, and a natural
 2-morphism $\alpha_t \: \ev_t\circ c_{\X} \Rightarrow \id_{\X}$.

Similarly,  we define  the {\bf free loop stack} of $\X$ to be $L\X:=\Map(S^1,\X)$. 
If $\X$ is topological (resp., paratopological), then so is $L\X$. There is, however, another
way of defining the free loop stack. Namely,  we can define $L\X$ to be
     $$\X\times_{\Delta,\X\times\X,(\ev_0,\ev_1)}P\X.$$
In contrast to the case of topological spaces,  these definitions are not expected to be 
equivalent in general. However, if we assume that $\X$ is a Hurewicz topological stack, 
then all ``reasonable'' definitions are equivalent. This is true thanks to Proposition \ref{P:glue}.

Assume now that $x \in \X$ is a basepoint. The {\bf based loop stack}
$\Omega_x\X$ of $\X$ is defined by 
  $$\Omega_x\X := * \times_{x,\X,\ev}L\X,$$
 where $\ev \: L\X \to \X$ is evaluation at the basepoint of $S^1$. Again,  
 as in the previous paragraph, if we do not assume that $\X$ is 
 Hurewicz, there may be more than one way of defining $\Omega_x\X$, and these
 definitions may not be equivalent. But if $\X$ is Hurewicz, all reasonable definitions
 will be equivalent.

\subsection{Homotopy fiber of a morphism of stacks}{\label{SS:htpyfiber}}

In classical topology there is a standard procedure for replacing an arbitrary 
continuous map $f \: X \to Y$ of topological spaces by a fibration. This construction 
involves taking path spaces and fiber products, both of which are available to us in the 
2-category of (para)topological stacks, thanks to Theorems \ref{T:mapping1}, 
\ref{T:mapping2} and Proposition \ref{P:arblimits}.

Let $f \: \X \to \Y$ be a morphism of stacks. 
Set $\tX:=\X\times_{f,\Y,\ev_1}P\Y$, where $\ev_1 \: P\Y \to \Y$
is the time $t=1$ evaluation map. Note that if $\X$ and $\Y$ 
are (para)topological stacks, then so is $\tX$. 
We define $p_f \: \tX \to \Y$ to
be the composition $\ev_0\circ\pr_2$, and $i_f \: \X \to \tX$ to be
the map whose first and second components are $\id_{\X}$ and
$c_{\Y}\circ f$, respectively.

The proof of the following proposition, as well as a thorough discussion of the notion 
of fibration of stacks, will appear elsewhere.

\begin{prop}{\label{P:replacement}}
 Notation being as above, we have a factorization $f=p_f\circ i_f$ such that:
 \begin{itemize}
   \item[$\mathbf{i.}$] the map $p_f \: \tX \to \Y$ is a Hurewicz
     fibration;

   \item[$\mathbf{ii.}$] the map $i_f \: \X \to \tX$ is a closed embedding
   and it admits a strong deformation retraction $r_f \: \tX \to \X$.
 \end{itemize}
\end{prop}

We remark that the above proposition is quite formal and is valid for all stacks. 
What is interesting, and that is where we make use of the results of this papers,
is that in the case where $\X$ and $\Y$ are (para)topological   
stacks, $\tX$ is again a (para)topological stack. This is 
important because it means that when working with (para)topological stacks 
the fibration replacement of morphisms keep us in  
the 2-category of (para)topological stacks.
 
We define the {\bf homotopy fiber} over a point $y \in \Y$  of a morphism $f \: \X \to \Y$ 
of (para)topological stacks to be
  $$\hFib_y(f):= *\times_{y,\Y,p_f}\tX,$$
where $p_f$ is as in Proposition \ref{P:replacement}. 
Since the 2-category of (para)topological stacks is closed under fiber products, 
the homotopy fiber of a morphisms of (para)topological stacks  is again
a (para)topological  stack.
 
\section{Homotopy invariance of mapping stacks}{\label{S:Invariance}}

In this section we study homotopy invariance of the mapping 
stack $\Map(\Y,\X)$  under change of $\X$. Although we can not 
make a statement in full generality, our result 
(Theorem \ref{T:invariance})  applies to many important classes 
of examples (e.g., loop stacks). A useful consequence of our invariance
theorem is that to calculate the homotopy type of the mapping stack 
$\Map(Y,\X)$, where $Y$ is a paracompact topological space and $\X$ is a 
paratopological stack, one can replace $\X$ with its classifying space $X$; 
see Corollary \ref{C:invariance1} below for the precise statement.

\subsection{Weakly parashrinkable morphisms}{\label{SS:weaklypara}}

There is a class of representable morphisms of stacks which shares most of the
nice properties of parashrinkable morphisms but is particularly better-behaved. 
Although not strictly necessary, we deemed appropriate to formulate our invariance 
theorem in terms of these morphisms.

\begin{defn}{\label{D:weaklypara}}
  We say that a representable morphism $f \: \X \to \Y$ of stacks is {\bf weakly 
  parashrinkable} if for every paracompact topological space $T$, and 
  every morphism $p \:  T \to \Y$, the space $\Map_p(T,\X)$ of lifts of $p$ to $\X$ 
  is weakly contractible and non-empty.
\end{defn}

By definition, $\Map_p(T,\X)$ is the fiber of the morphism of 
stacks $f_* \: \Map(T,\X) \to \Map(T,\Y)$ over the point in  $\Map(T,\Y)$ corresponding to 
$p$. Note that $\Map_p(T,\X)$ is 
naturally equivalent 
to the space of sections of the projection map $T\times_{\Y}\X \to T$. In particular,
$\Map_p(T,\X)$ is an honest topological space.

\begin{lem}{\label{L:weaklypara}}
  The following statements are true about weakly parashrinkable morphisms:
  \begin{itemize}
  
    \item[$\mathbf{1.}$] Every parashrinkable morphism is weakly parashrinkable. 
    \item[$\mathbf{2.}$] Every weakly parashrinkable morphism is a universal weak 
         equivalence.
    \item[$\mathbf{3.}$] Base extension of a weakly parashrinkable morphism is
        weakly parashrinkable.      
  \end{itemize}

\end{lem}

\begin{proof}  Proof of ($\mathbf{1}$) is easy because the space of sections of the shrinkable 
    morphism $T\times_{\Y}\X \to T$ is clearly contractible. Parts ($\mathbf{2}$) and 
    ($\mathbf{3}$) are straightforward. 
\end{proof}

\begin{rem}{\label{R:pseudo}}
  The above discussion can be repeated for {\em weakly shrinkable} morphisms and
  {\em weakly pseudoshrinkable} morphisms  (\cite{Homotopytype},  
  Definition 5.1) as well. Since we do not need these concepts in this paper, we will
  avoid further discussion.
\end{rem}

\subsection{Invariance theorem}{\label{SS:invariance}}

\begin{thm}{\label{T:invariance}}
  Let $f \: \X' \to \X$ be a representable morphism of stacks. Assume that $f$ is 
  weakly parashrinkable (Definition  \ref{D:weaklypara} ). Let $Y$ be a paracompact 
  topological space. Then $f_*\:  \Map(Y,\X') \to \Map(Y,\X)$ is a universal 
  weak equivalence. If $Y$ is compact, then $f_*$ is weakly parashrinkable.
\end{thm}

\begin{proof}
    Let $T$ be a finite CW complex (or any compact topological space). 
    We show that for every morphism 
    $p \: T \to \Map(Y,\X)$, the space $\Map_p(T,\Map(Y,\X'))$ of lifts of $p$ to
    $\Map(Y,\X')$ is weakly contractible. It is easy to see that this implies that 
    $f_*\:  \Map(Y,\X') \to \Map(Y,\X)$ is a universal weak equivalence. 
     
     To show that $\Map_p(T,\Map(Y,\X'))$ is weakly contractible, observe that
     we have a natural homeomorphism
       $$\Map_p(T,\Map(Y,\X')) \cong \Map_{\tilde{p}}(T\times Y, \X'),$$
    where  $\tilde{p} \: T\times Y \to \X$ is the map corresponding to
    $p$ under the exponential law (Lemma \ref{L:exponential}). Since $T\times Y$,
    being the product a compact space and a paracompact space, is paracompact,
    and since $\X' \to \X$
    is weakly parashrinkable,  the right hand side of the equation is contractible.     
    
    In the case where $Y$ is compact, one can repeat the same argument as above,
    with $T$ an arbitrary paracompact topological space. One finds that 
    $f_*\:  \Map(Y,\X') \to \Map(Y,\X)$ is weakly parashrinkable.
\end{proof}

 
\begin{cor}{\label{C:invariance1}}
   Let $Y$ be a paracompact topological space and $\X$ a paratopological stack. Let
   $X$ be a classifying space for $\X$ with $\varphi \: X \to \X$ weakly parashrinkable
   (such an $X$ always exists, see Proposition \ref{P:classifying}). Then, 
   $\Map(Y,X)$ is a classifying space for $\Map(Y,\X)$, hence it represents the weak
   homotopy type of $\Map(Y,\X)$.   
\end{cor}

\begin{cor}{\label{C:invariance2}}
   Let $\X$ be a paratopological stack and $X$ a classifying space for it, 
   with $\varphi \: X \to \X$ weakly parashrinkable
   (such an $X$ always exists, see Proposition \ref{P:classifying}). Then, the free 
   loop space $LX$ of $X$ is a classifying space for the free loop space $L\X$ of $\X$ (hence, 
   has the  same weak homotopy type). In fact, we can arrange for the map $LX \to L\X$ to be  
   weakly parashrinkable (hence, a universal weak equivalence).
\end{cor}

Of course, there is nothing special about $S^1$ in the above corollary. One can take any compact 
topological space instead of $S^1$.

\begin{cor}{\label{C:composition}}
    Composition of two weakly parashrinkable morphisms is   weakly parashrinkable.
\end{cor}

\begin{proof}
  Let $f \: \X \to \Y$ and $g \: \Y \to \Z$ be weakly parashrinkable morphisms. Let $T$ be a 
  paracompact topological space and $p \: T \to \Z$ a morphism. We want to show that
  $\Map_p(T,\X)$ is weakly contractible. We have a cartesian square
        $$\xymatrix{ \Map_p(T,\X) \ar[r] \ar[d]_{f_*} & \Map(T,\X) \ar[d]^{f_*} \\
           \Map_p(T,\Y)  \ar[r] & \Map(T,\Y)  }$$
  By Theorem \ref{T:invariance}, the right vertical arrow is a universal weak 
  equivalence.  Hence, so is the left vertical arrow. 
  Since $g$ is weakly parashrinkable, $ \Map_p(T,\Y)$
  is weakly contractible. Therefore,   $\Map_p(T,\X)$ is weakly contractible.
\end{proof}

The above corollary reveals a main advantage of weakly parashrinkable morphisms over
parashrinkable morphisms: they are closed under composition. The same thing is 
true with weakly shrinkable morphisms and weakly pseudoshrinkable morphisms 
as well (see Remark \ref{R:pseudo}).

  \begin{rem}
  Homotopy invariance of $\Map(\Y,\X)$ with respect to $\Y$ is not an interesting problem. 
  To illustrate this, consider the case where $\Y=\mathcal{B}G$ is the classifying stack 
  of a group $G$, and $\Y'=|BG|$ is the Milnor classifying space of $G$. Let $\X=X$ be an 
  arbitrary topological space. There is a natural parashrinkable morphism $BG \to \mathcal{B}G$.
  The induced map
       $$\Map(\mathcal{B}G,X) \to \Map(BG,X)$$ 
  is in general far from being
  a weak equivalence. The left hand side is nothing but $\Map(*,X)=X$, while the 
  right hand side could be considerably more complicated.
\end{rem}

\section{An example: loop stack of $\mathcal{B}G$ and twisted loop groups}{\label{S:Loop}}

In this section,  we describe the loop stack $L\mathcal{B}G$ of the classifying stack
of a topological group $G$. We show that the homotopy type of $L\mathcal{B}G$
can be calculated in terms 
of the twisted loop groups of $G$. In the case where $G$ is discrete,  
$L\mathcal{B}G$ is equivalent to the inertia stack of $\mathcal{B}G$. See 
Theorem \ref{T:loopstack}.

In what follows, all $G$-torsors are right $G$-torsors.

\begin{lem}{\label{L:Milnor}}
  Let $G$ be a topological group and $BG$ its Milnor classifying space. Let $T$ be a 
  paracompact topological space. Then, every $G$-torsor on $[0,1]\times T$ is isomorphic
  to the pull-back of a $G$-torsor on $T$.
\end{lem}

\begin{proof}
  It is well-known that isomorphism classes of
  $G$-torsors  on a paracompact space $T$ are in natural bijection with 
  homotopy classes of  maps $T \to BG$. Apply this fact to $T$ and 
  $[0,1]\times T$, both of which are paracompact, and
  use the fact that the projection $[0,1]\times T \to T$ is a homotopy equivalence.
\end{proof}

\begin{lem}{\label{L:gluetorsor}}
  Let $U$ be a topological space and $f \: U \to G$ a continuous map. Let $T_f$ 
  be the quotient space of   $G\times U\times [0,1] $ obtained by gluing 
  $G\times U\times\{0\}$ to $G\times U\times\{1\}$ along the map 
  $(g,u,0) \mapsto \big(f(u)g,u,1\big)$. Then, $T_f$ is naturally a (right) $G$-torsor over
  $U\times S^1$. Furthermore, if $f' \: U \to G$ is another continuous map, then 
  $T_{f'}$ is isomorphic to $T_f$ as a $G$-torsor if and only if there is a function 
  $\delta \: U \to G$ such that $f'$ is homotopic to the conjugate $\delta f \delta^{-1}$ of
  $f$ under $\delta$.
\end{lem}

\begin{proof}
   A straightforward argument shows that a morphism $T_f \to T_{f'}$  of $G$-torsors 
   corresponds precisely to a map $\gamma \: U\times [0,1] \to G$ such that 
       $$f'(u)\gamma(u,0)=\gamma(u,1)f(u).$$
   If such a $\gamma$ exists, then $h(u,t):=f'(u)\gamma(u,t)\gamma(u,1)^{-1}$ gives
   a homotopy  from $\delta f\delta^{-1}$ to $f'$, where $\delta(u):=\gamma(u,1)$. 
   Conversely, given a homotopy $h(u,t)$ from $\delta f\delta^{-1}$ to $f'$, for some
   $\delta \: U \to G$, we define $\gamma(u,t):=f'(u)^{-1}h(u,t)\delta(u)$. 
\end{proof}

Let $\pi_0(G)$ denote the group of path components of $G$, and $C_G$ the 
set of conjugacy classes of elements of $\pi_0(G)$. Choose a set of 
representative $\alpha_i \in G$, $i\in C_G$, of conjugacy classes 
of path components of $G$.   
Let $T_{\alpha} $ be the $G$-torsor over $S^1$ obtained from the trivial
 torsor $G\times [0,1] $ on $[0,1]$ by gluing
$G\times\{0\}$ to $G\times\{1\}$ along the map $(g,0) \mapsto (\alpha g,1)$,
as in Lemma \ref{L:gluetorsor}. 
The automorphism group of $T_{\alpha}$ is isomorphic, as a topological group, 
to the  {\bf twisted loop group}   $L_{(\alpha)}G$ (\cite{PrSe}, $\S$3.7)  
associated to the conjugation action
$g \mapsto \alpha g \alpha^{-1}$ of $\alpha$ on  
$G$.  In other words,
       $$L_{(\alpha)}G:=\{\gamma \: \mathbb{R} \to G \ | \ 
                           \gamma(\theta+1)=\alpha\gamma(\theta)\alpha^{-1}\}.$$ 
(Note that our setting is slightly more general than that of [ibid.] because we are not
assuming that $G$ is connected.)                           
In the case where $\alpha$ belongs to the path component of the identity,   
$L_{(\alpha)}G$ is isomorphic to the loop group $LG=\Map(S^1,G)$. 
When $G$ is discrete, $L_{(\alpha)}G$ is the centralizer of $\alpha$.

For every $\alpha$ as above, let $p_{\alpha} \: * \to L\mathcal{B}G$ denote the 
map corresponding to the $G$-torsor $T_{\alpha}$. We have the following.

\begin{lem}{\label{L:product}}
  Let $p_{\alpha} \: * \to L\mathcal{B}G$ be defined as above. 
  Then, we have a natural
  isomorphism of topological groups
       $$*\times_{p_{\alpha},L\mathcal{B}G,p_{\alpha}}*\cong    L_{(\alpha)}G.$$
  If $\alpha, \beta \in G$ map to different elements in $C_G$, then 
        $$*\times_{p_{\alpha},L\mathcal{B}G,p_{\beta}}*\cong \emptyset.$$
 \end{lem}

 \begin{proof}
   The first statement follows from fact that the automorphism group of the 
   $G$-torsor $T_{\alpha}$ is isomorphic to $L_{(\alpha)}G$ as a topological group. 
   The second statement   follows from the fact that,  $T_{\alpha} $ and 
   $T_{\beta} $ are isomorphic as $G$-torsors  if and only if the images
   of $\alpha$ and $\beta$ in $\pi_0(G)$ belong to the same conjugacy class
   (see Lemma \ref{L:gluetorsor}).
 \end{proof}

\begin{lem}{\label{L:torsor}}
  The map 
     $P=\coprod_{i\in C_G} p_{\alpha_i} \:  \coprod_{i\in C_G} * \to L\mathcal{B}G$, 
  where $*$ stands for a point, has the property ($\mathbf{A2}$) of 
  Definition \ref{D:paratop}  for every locally contractible  paracompact
   topological space $W$ (in particular, for every finite CW complex $W$).
\end{lem}

\begin{proof}
  Pick a map $h\: W \to L\mathcal{B}G$, and let $T \to W\times S^1$ be 
  the corresponding  $G$-torsor. We want to show that every point $x \in W$ 
  has an open neighborhood $U$ such that $h|_U$ lifts along some 
  $p_{\alpha_i} \: * \to L\mathcal{B}G$,  that is, 
  $T|_{U\times S^1}$ is isomorphic to the $G$-torsor 
  $U\times T_{\alpha_i}$, for some $i \in C_G$. We claim that any contractible 
  neighborhood  $U$ of $x$ has the desired property. By Lemma \ref{L:Milnor}, 
  we know that $T$  is of the form $T_f$ for some $f \: W \to G$ (see Lemma 
  \ref{L:gluetorsor} for notation). Since $U$ is contractible, $f|_U$ is homotopic 
  to a constant map into some path component of $G$.
  Hence, $f|_U$ is homotopic to a conjugate of the constant map 
  $\alpha_i \: U \to G$,  for some $\alpha_i$. Therefore, by Lemma 
  \ref{L:gluetorsor}, $T|_U=T_{f|_U}$ is isomorphic
  to $U\times T_{\alpha_i}$ as a $G$-torsor. In other words, $f|_U$ lifts along  
  $p_{\alpha_i} \: * \to L\mathcal{B}G$.
\end{proof}

\begin{cor}{\label{C:loopstack}}
     The map $P$ of Lemma \ref{L:torsor} induces a natural map 
       $$\bar{P}\: \coprod_{i\in C_G} \mathcal{B}L_{(\alpha_i)}G 
                                  \to  L\mathcal{B}G$$
     which has the property of Lemma \ref{L:approximation}  for every locally 
     contractible paracompact topological space $W$ 
     (in particular, for every finite CW complex $W$). In particular, $\bar{P}$ is a 
     universal weak equivalence.
\end{cor}

\begin{proof}
  For the existence of the map $\bar{P}$ and that it has the property
   of Lemma \ref{L:approximation} see the proof of Lemma \ref{L:approximation} 
  (and the remark after the proof). Observe that we are also using 
  Lemma \ref{L:product}. 
\end{proof}

Finally, we put together what we have proved about the loop stack of $\mathcal{B}G$
in the following.

\begin{thm}{\label{T:loopstack}}
   Let $G$ be a topological group. Then, there are natural weak homotopy equivalences
     $$LBG\cong L\mathcal{B}G \cong \coprod_{i\in C_G} BL_{(\alpha_i)}G.$$ 
   In particular, when $G$ is discrete, there are natural weak homotopy equivalences
     $$LBG\cong \mathcal{IB}G \cong G\times_GEG,$$ 
   where  $\mathcal{IB}G$ is the inertia stack and $G\times_GEG$ is
   the Borel construction applied to the conjugation action of $G$ on itself.
   (The equivalence $LBG\cong \mathcal{IB}G$ 
   is indeed an equivalence of stacks. This can be easily verified directly.)
\end{thm}

We can  describe the weak homotopy type of the based loop stack 
$\Omega\mathcal{B}G$ in a similar fashion. (Observe that $\mathcal{B}G$ is
a Hurewicz topological stack, so by the discussion of $\S$ \ref{SS:path} the
different definitions for the based loop stack agree.) Consider the  map
$\ev \:  L\mathcal{B}G \to \mathcal{B}G$ which evaluates a loop at its basepoint.
The based loop stack $\Omega\mathcal{B}G$ is the fiber of $\ev$ over the base 
point of $\mathcal{B}G$ corresponding to the trivial $G$-torsor. Note that the 
composite map 
  $$\ev\circ \bar{P} \: \coprod_{i\in C_G} \mathcal{B}L_{(\alpha_i)}G 
                                  \to \mathcal{B}G$$
is the one induced by the evaluation maps $\ev_0 \: L_{(\alpha_i)}G \to G$,
$\gamma \mapsto \gamma(0)$. Therefore, the fiber of $\ev\circ \bar{P}$
over the basepoint of $\mathcal{B}G$ is
           $$\coprod_{i\in C_G} \mathcal{B}\Omega_{(\alpha_i)}G,$$
where   $\Omega_{(\alpha)}G$ is the subgroup of the twisted loop group
$L_{(\alpha)}G$ defined by
      $$\Omega_{(\alpha)}G:=\{\gamma \: \mathbb{R} \to G \ | \ 
             \gamma(\theta+1)=\alpha\gamma(\theta)\alpha^{-1}, \ \gamma(0)=1_G\}.$$  
Note that $\Omega_{(\alpha)}G$ only sees the path component of the identity.
In other words, $\Omega_{(\alpha)}G=\Omega_{(\alpha)}(G^0)$. In particular, if
$G$ is discrete, then $\Omega_{(\alpha)}G$ is trivial.

We have the following 2-cartesian diagram
  $$\xymatrix{  \coprod_{i\in C_G} \mathcal{B}\Omega_{(\alpha_i)}G  \ar[d]_Q \ar[r] 
                         &  \coprod_{i\in C_G} \mathcal{B}L_{(\alpha_i)}G \ar[d]^{\bar{P}} \\
                               \Omega\mathcal{B}G \ar[r]              &      L\mathcal{B}G}$$
                                  
The morphism $Q$, being the base extension of $\bar{P}$ inherits the same property 
that $\bar{P}$ does in Corollary \ref{C:loopstack}. In particular, $Q$ is a universal
weak equivalence. We have the following.

\begin{thm}{\label{T:pointed}}
     Let $G$ be a topological group. Then, there are natural weak homotopy equivalences
     $$\Omega BG\cong \Omega\mathcal{B}G \cong 
                          \coprod_{i\in C_G} B\Omega_{(\alpha_i)}G.$$ 
   In particular, when $G$ is discrete, $\Omega \mathcal{B}G$ is weakly 
   homotopy equivalent to  the discrete set $C_G$. (The latter is indeed an equivalence
   of stacks as can be easily verified directly.)
\end{thm}

\begin{proof}
    Everything follows immediately from the discussion  above
    (thanks to Proposition \ref{P:Haefliger}),
    except for the weak equivalence $\Omega BG\cong \Omega\mathcal{B}G$.
    To prove this, consider the 2-commutative diagram
      $$\xymatrix{ LBG \ar[r]^{\varphi_*} \ar[d]_{\ev} & 
                                       L\mathcal{B}G \ar[d]^{\ev} \\
                          BG \ar[r]_{\varphi}  & \mathcal{B}G}$$
      Let $V \subset BG$ be the fiber of     $\varphi \: BG \to    \mathcal{B}G$
      over the base point of $\mathcal{B}G$ and $W=\ev^{-1}(V)$ the fiber of
      $\varphi \circ \ev \:  LBG \to \mathcal{B}G$.  Observe that $V$ is contractible
      (because $\varphi$ is parashrinkable), and that $\ev    \:  LBG \to BG$
      is a fibration. Therefore, the inclusion $j \: \Omega BG \hookrightarrow W$ is a 
      weak homotopy equivalence. On the other hand, by the 2-commutativity of the
      above diagram, $W$ fits in the 2-cartesian diagram
              $$\xymatrix{W \ar[r]^f \ar[d] & \Omega \mathcal{B}G \ar[d] \\
                          LBG \ar[r]_{\varphi_*} & 
                                       L\mathcal{B}G} $$
       By      Corollary \ref{C:invariance2} and Proposition \ref{P:Haefliger}, $f$ is 
       a  weak equivalence.
       Therefore, the composition    $f\circ j \:   \Omega BG \to   \Omega \mathcal{B}G$
       is also a weak equivalence.   This completes the proof of the theorem.
\end{proof}

One can use Theorem \ref{T:pointed} and apply the loop space/stack functor 
repeatedly to find a description  of the week homotopy types of the mapping 
stacks $\Map(S^n,\mathcal{B}G)$ and also the mapping spaces  $\Map(S^n,BG)$ 
in terms of certain twisted mapping groups of $G$. We will leave it to the interested 
reader to work out the details.

\providecommand{\bysame}{\leavevmode\hbox
to3em{\hrulefill}\thinspace}
\providecommand{\MR}{\relax\ifhmode\unskip\space\fi MR }
\providecommand{\MRhref}[2]{%
  \href{http://www.ams.org/mathscinet-getitem?mr=#1}{#2}
} \providecommand{\href}[2]{#2}


\begin{thebibliography}{10}



\bibitem[BGNX]{BGNX} K.~Behrend, G.~Ginot, B.~Noohi, P.~Xu, \emph{String topology for stacks},
    arXiv:0712.3857v1 [math.AT].

\bibitem[Ch]{Chen} W.~Chen, \emph{On a notion of maps between orbifolds, I. Function
spaces},  Commun. Contemp. Math.  \textbf{8}  (2006),  no. 5,
569--620.

\bibitem[EbGi]{EbGi} J.~Ebert, J.~Giansiracusa, \emph{Pontrjagin-Thom maps and the homology of the moduli stack
of stable curves}, arXiv:0712.0702v2 [math.AT].

\bibitem[LuUr]{LuUr} E.~Lupercio, B.~Uribe, \emph{Loop groupoids, gerbes, and twisted sectors on orbifolds},
Orbifolds in mathematics and physics (Madison, WI, 2001),  163--184, Contemp. Math., 310, Amer. Math. Soc.,
Providence, RI, 2002.

\bibitem[No1]{Foundations}  B.~Noohi, \emph{Foundations of topological
 stacks, I}, math.AG/0503247v1.

\bibitem[No2]{Homotopytype}  B.~Noohi, \emph{Homotopy types of topological
stacks}, arXiv:0808.3799v1 [math.AT].

\bibitem[PrSe]{PrSe}  A.~Pressley, G.~ Segal, \emph{Loop groups}, Oxford Mathematical Monographs, New York, 1986.

\bibitem[Wh]{Whitehead}
 G.~W.~Whitehead, \emph{Elements of Homotopy Theory}, Graduate Texts
 in Mathematics, 61. Springer-Verlag, New York-Berlin, 1978.

\end{thebibliography}
\end{document}